\documentclass[11pt]{article}
\usepackage{amsfonts}
\usepackage{amssymb}
\usepackage{amssymb,amsfonts,amsmath,amsthm,cite, color}
\usepackage{epsfig}
\newtheorem{thm}{Theorem}[section]

\newtheorem{lem}{Lemma}[section]
\newtheorem{prop}{Proposition}[section]
\newtheorem{defn}{Definition}[section]

\makeatletter \@addtoreset{equation}{section}

\def\pf{\noindent {\it Proof.\ }}
\def\qed{\hfill \rule{4pt}{7pt}}

\begin{document}

\begin{center}
{{\Large\bf Congruences for the Number of Cubic Partitions

 Derived from Modular Forms}}

\vskip 6mm

{ William Y.C. Chen$^1$ and Bernard L.S. Lin$^2$
\\[%
2mm] Center for Combinatorics, LPMC-TJKLC\\
Nankai University, Tianjin 300071,
P.R. China \\[3mm]
$^1$chen@nankai.edu.cn, $^2$linlishuang@cfc.nankai.edu.cn \\[0pt%
] }
\end{center}

\begin{abstract}
We obtain congruences for the number $a(n)$ of cubic partitions
using modular forms. The notion of cubic partitions is introduced by
Chan and named by Kim in connection with Ramanujan's cubic continued fractions.
 Chan has shown that $a(n)$ has several analogous properties to
 the number $p(n)$ of partitions, including the generating function,
 the continued fraction, and congruence relations.
To be more specific, we show that  $a(25n+22)\equiv 0\ ({\rm mod}\
5)$, $a(49n+15) \equiv a(49n+29) \equiv a(49n+36) \equiv a(49n+43)
\equiv 0 \ ({\rm mod}\ 7)$. Furthermore, we prove that $a(n)$ takes
infinitely many even values and infinitely odd values.
\end{abstract}



\section{Introduction}

The main objective of this paper is to study congruence relations
for the number of cubic partitions by constructing suitable
modular forms. The number of cubic partitions, denoted by $a(n)$,
originated from the work of Chan \cite{Chan08a} in connection
with Ramanujan's cubic continued fraction which is often denoted by
\[
G(q):=\frac{q^{1/3}}{1}_+\frac{q+q^2}{1}_+\frac{q^2+q^4}{1}_+\frac{q^3+q^6}{1}_{+\cdots}
,\quad |q|<1.
\]
On page $366$ of his Lost Notebook, Ramanujan
claimed that there are many properties of $G(q)$ which are
analogous to  Rogers-Ramanujan continued fraction $R(q)$ \cite{Ramanujan88}
\[
R(q):=\frac{q^{1/5}}{1}_+\frac{q}{1}_+\frac{q^2}{1}_+\frac{q^3}{1}_{+\cdots},
\quad |q|<1.
\]

Motivated by Ramanujan's observation,
many new results on $G(q)$ analogous to those for $R(q)$ have been found,
see, e.g.,  Chan \cite{HChan}.
 To give an overview of recent results on $a(n)$, it is informative to
  recall relevant background on the generating function of $p(n)$
 and the  Rogers-Ramanujan continued fraction $R(q)$.

Ramanujan obtained many theorems on $R(q)$, see Andrews and
Berndt \cite{Andrews}. In particular, he  discovered the following
beautiful identities on  $R(q)$ and $1/R(q)$.
\begin{eqnarray}
\frac{1}{R(q)}-1-R(q)&=&\frac{(q^{1/5};q^{1/5})_\infty}{q^{1/5}(q^5;q^5)_\infty}\label{R1}\\[5pt]
\frac{1}{R^5(q)}-11-R^5(q)&=&\frac{(q;q)_\infty^6}{q(q^5;q^5)_\infty^6}.\label{R2}
\end{eqnarray}
Here $(q;q)_\infty$ is the usual notation for
$\prod\limits_{n=1}^\infty (1-q^n)$.

Berndt \cite[p.165]{Berndt04} gave a beautiful proof of the
following classical identity of Ramanujan by using the continued
fraction $R(q)$:
\begin{equation}\label{rami}
\frac{(q;q)_\infty ^6}{(q^5;q^5)_\infty ^5}\sum_{n=0}^\infty
p(5n+4)q^n=5.
\end{equation}
 Dividing  \eqref{R2}  by
\eqref{R1}, we get
\begin{eqnarray}\lefteqn{
\frac{(q;q)_\infty ^6}{q^{4/5}(q^{1/5};q^{1/5})_\infty
(q^5;q^5)_\infty ^5}} \nonumber \\[6pt]  \label{bd}
& =R^4(q)-R^3(q)+2R^2(q)-3R(q)+5+\frac{3}{R(q)}
+\frac{2}{R^2(q)}+\frac{1}{R^3(q)}+\frac{1}{R^4(q)}.
\end{eqnarray}
Now, (\ref{rami}) can be easily deduced from (\ref{bd}) by
extracting the integer powers of $q^n,n\geq 0$ from both sides of
above identity since $R(q)$ has only terms in the form of
$q^{n+1/5}$. Ramanujan's congruence on $p(n)$ modulo $5$ can be
derived directly from \eqref{rami}
\begin{equation}\label{dyson}
p(5n+4)\equiv 0\ ({\rm mod}\ 5).
\end{equation}

Recently, using two identities of Ramanujan \cite{Ramanujan88} on $G(q)$, see
also Berndt \cite[p.345, Entry 1]{Berndt91}, Chan \cite{Chan08a}
has found the following  identities on $G(q)$ and $1/G(q)$
analogous to the above  identities \eqref{R1} and \eqref{R2}:
\begin{eqnarray}
\frac{1}{G(q)}-1-2G(q)&=&\frac{(q^{1/3};q^{1/3})_\infty
(q^{2/3};q^{2/3})_\infty}
{q^{1/3}(q^3;q^3)_\infty(q^6;q^6)_\infty} , \label{G1}\\[5pt]
\frac{1}{G^3(q)}-7-8G^3(q)&=&\frac{(q;q)_\infty ^4(q^2;q^2)_\infty
^4}{q(q^3;q^3)_\infty ^4(q^6;q^6)_\infty ^4}.\label{G2}
\end{eqnarray}

 Motivated by the idea of Berndt, Chan derived the following identity
 by dividing both sides of \eqref{G2} by
 \eqref{G1} and then setting $q\rightarrow q^3$:
 \begin{eqnarray}
 \frac{1}{(q;q)_\infty(q^2;q^2)_\infty}=q^2\frac{(q^9;q^9)_\infty^3
 (q^{18};q^{18})_\infty^3}{(q^3;q^3)_\infty^4(q^6;q^6)_\infty^4}\left(4G^2(q^3)-2G(q^3)
 +3+\frac{1}{G(q^3)}+\frac{1}{G^2(q^3)}\right). \label{ChanGDivde}
 \end{eqnarray}

Observing that  the powers of $q$ in $G(q^3)$ are in the form of $3n+1$, we find
 \[
 \sum_{n=0}^\infty \big[q^{3n}\big]\left(4G^2(q^3)-2G(q^3)
 +3+\frac{1}{G(q^3)}+\frac{1}{G^2(q^3)}\right) q^{3n} = 3.
 \]
 It is now natural to define a function $a(n)$  by the left hand side of
\eqref{ChanGDivde}
\begin{equation}\label{Defa}
\sum_{n=0}^\infty a(n)q^n=\frac{1}{(q;q)_\infty(q^2;q^2)_\infty},
\end{equation}
 and it is natural to expect $a(n)$  to have analogous properties to $p(n)$.

 Extracting those terms whose powers of $q$ are in the form of $3n+2$
on both sides of \eqref{ChanGDivde}, and then simplifying and
setting $q^3\rightarrow q$, Chan established the following elegant
identity analogous to \eqref{rami}
\begin{equation}\label{ChanIdena}
\sum_{n=0}^\infty a(3n+2)q^n=3\frac{(q^3;q^3)_\infty^3
(q^6;q^6)_\infty ^3}{(q;q)_\infty ^4 (q^2;q^2)_\infty ^4}.
\end{equation}
The above identity immediately leads to the following congruence
\begin{equation}\label{Conga}
a(3n+2)\equiv 0\ ({\rm mod}\ 3),
\end{equation}
which is analogous to Ramanujan's congruence \eqref{dyson} for
$p(n)$.

From the point of view of partitions, it is obvious from the
generating function \eqref{Defa} that $a(n)$ is the number of
partition pairs $(\lambda,\mu)$ where $|\lambda|+|\mu|=n$ and $\mu$
only has even parts.  Chan has called $a(n)$ a certain partition
function.  Kim \cite{Kim08} called  such partitions counted by
$a(n)$ cubic partitions owing to the fact that $a(n)$ is
close related to Ramanujan's cubic continued fraction.

Based on the cubic partition interpretation of $a(n)$, Chan \cite{Chan08c} asked whether
there exist a function analogous to Dyson's rank  that leads to
 a combinatorial interpretation  of the congruence \eqref{dyson}.
Kim  \cite{Kim08} discovered  a crank function
$N_V^{a}(m,n)$ for cubic partitions. Let $M'(m,N,n)$ be the number of cubic
partitions of $n$ with crank $\equiv m \  ({\rm mod}\ N)$, Kim
proved that
\[
M'(0,3,3n+2)\equiv M'(1,3,3n+2) \equiv M'(2,3,3n+2)\  ({\rm mod}\
3),
\]
which implies  \eqref{Conga}.

Our main results are concerned with congruences for $a(n)$ modulo $5$ and $7$
which are in the spirit of Ramanujan's classical congruences modulo $7$ and $11$.
Recall that Ramanujan obtained more general  congruences modulo $5^k$:
\begin{equation}\label{gencongp}
 p(5^kn+r_k)\equiv
0\ ({\rm mod}\ 5^k),
\end{equation}
where $k\geq 1$ and
$5^k r_k \equiv 1\ ({\rm mod}\ 24)$. In analogy with Ramanujan's congruences,
 Chan  considered the general
congruences for $a(n)$ modulo powers of $3$. Employing the
method of Hirschhorn and Hunt \cite{Hirschhorn} to
prove \eqref{gencongp},  Chan \cite{Chan08b} derived the following
congruence as a consequence of \eqref{ChanIdena}.

\begin{thm}
For $k\geq 1$,
\begin{equation}\label{genconga}
a(3^kn+c_k)\equiv  0\ ({\rm mod}\ 3^{k+\delta(k)}),
\end{equation}
where $c_k$ is the reciprocal modulo $3^k$ of $8$, and $\delta(k)=1$
if $k$ is even and $\delta(k)=0$ otherwise.
\end{thm}

In the general case, Ramanujan conjectured that there are only three
choices for a prime $l$ such that the congruence  $p(ln+c) \equiv 0\
({\rm mod}\ l)$ holds,  namely, $l=5, 7, 11$. This conjecture has
been confirmed  by Ahlgren and Boylan \cite{Alhgren03} based on the
work of Kiming and Olsson \cite{Kiming}.  Chan \cite{Chan08b} raised
the problem of finding simple congruences for $a(n)$ besides
$a(3n+2)\equiv 0\mod 3$.  Recently, Sinick \cite{Sinick08} has shown
that there does not exist
 other primes $l$ such that $a(ln+c)\equiv 0 \ ({\rm
mod}\ l)$ except that $l=3$. In analogy with the results for $p(n)$ due
to Ono \cite{Ono00} and Ahlgren \cite{Alhgren00}, Chan
\cite{Chan08c} obtained the following theorem concerning congruences
for $a(n)$ modulo powers of a prime.
\begin{thm}\label{thm1.2}
Let $m\geq 5$ be prime and $j$ a positive integer. Then a positive
proportion of the primes $Q \equiv -1 \ ({\rm mod}\ 128m^j)$ have
the property that
\[
a\left(\frac{mQn+1}{8}\right)\equiv 0 \ ({\rm mod}\ m^j),
\]
for every $n$ coprime to $Q$.
\end{thm}

The above theorem implies that   for every integer $n$ there exists
infinitely many non-nested arithmetic progressions $An+B$ for prime
$m\geq 5$ and positive integer $j$ such that
\[
a(An+B)\equiv  0\ ({\rm mod}\ m^j).
\]

It should be noted that although the proof of  Theorem \ref{thm1.2} leads to some Ramanujan-type congruences
modulo $m^j$, it does not cover all the congruences in form
of $a(An+B)\equiv  0\ ({\rm mod}\ m^j)$.
 Chan
\cite{Chan08b} studied the case for $m=3$, which is not in the scope
of Theorem \ref{thm1.2}. This paper is devoted to finding concrete
congruences for the cases $m=5,7$ and $j=1$, which are also out of
the range of Theorem \ref{thm1.2} since the  $Q$ is larger than
$1278$ and $2686$ for $m=5$ and $7$, respectively. To be precise, we
derive the following congruences by constructing suitable modular
forms.

\begin{thm}\label{thm1.3} For every nonnegative integer $n$, we have
\begin{eqnarray*}
a(25n+22)&\equiv& 0\ ({\rm mod}\ 5).
\end{eqnarray*}
\end{thm}

It would  be interesting to  give a combinatorial interpretation of
the above congruence by finding a suitable crank function. In the
following theorem, we present some congruences modulo $7$.

\begin{thm}\label{thm1.4}For every nonnegative integer $n$, we have
\begin{eqnarray*}
a(49n+15)\equiv a(49n+29)\equiv a(49n+36)\equiv a(49n+43)\equiv 0\
({\rm mod}\ 7).
\end{eqnarray*}
\end{thm}

The last section of this paper is focused on the parity of $a(n)$.
Recall that  Kolberg  \cite{Kolberg} has shown
 that $p(n)$ takes both even and odd values infinitely often.
From numerical evidence,   we conjecture
that when $n$ tends to infinity the parities of $a(1), a(2), \ldots, a(n)$ are
equidistributed. While we have not
been able to prove this conjecture, we shall show that there are
infinitely many even values of $a(n)$ and there are infinitely many
odd values of $a(n)$.

\section{Preliminaries}

To make this paper self-contained, we give an overview of the
background relevant to the proofs of the congruences for $a(n)$ by
using modular forms. For more details on the theory of modular
forms, see for example, Koblitz \cite{Koblitz} and Ono \cite{Ono04}.

For a rational integer $N\geq 1$, the congruence subgroup
$\Gamma_0(N)$ of $SL_2(\mathbb{Z})$ is defined  by
\[
\Gamma_0(N):=\left\{\left(
                 \begin{array}{cc}
                   a & b \\
                   c & d \\
                 \end{array}
               \right)\Big|c\equiv 0\ ({\rm mod}\ N)
\right\}.
\]
Let $\gamma=\left(
                 \begin{array}{cc}
                   a & b \\
                   c & d \\
                 \end{array}
               \right)\in SL_2(\mathbb{Z})$ act on the complex upper
half plane
\[
\mathbb{H}:=\{z\in\mathbb{C}|{\rm Im}(z)>0\}
\] by the linear
fractional transformation
\[
\gamma z:=\frac{az+b}{cz+d}.
\]
Suppose that $k$ is a positive integer and $\chi$ is a Dirichlet
character modulo $N$.

\begin{defn}
Let $f(z)$ be a holomorphic function on $\mathbb{H}$ and satisfy the
following relation for all $\gamma \in \Gamma_0(N)$ and all $z\in
\mathbb{H}$,
\[
f(\gamma z)=\chi(d)(cz+d)^kf(z).
\]
 In addition, if  $f(z)$ is also holomorphic at the cusps of
$\Gamma_0(N)$, we call such a function $f(z)$ a {\it modular form} of
weight $k$ on $\Gamma_0(N)$.
\end{defn}

 The modular forms of weight $k$ on
$\Gamma_0(N)$ with Dirichlet character $\chi$ form a
finite-dimensional complex vector space denoted by
$M_k(\Gamma_0(N),\chi)$. For convenience, we write $M_k(\Gamma_0(N))$
for $M_k(\Gamma_0(N),\chi)$ when $\chi$ is the trivial Dirichlet
character.

 Dedekind's eta function is defined by
\[
\eta(z):=q^{\frac{1}{24}}\prod_{n=1}^\infty (1-q^n),
\]
where $q=e^{2\pi iz}$ and ${\rm Im}(z)>0$. It is well-known that
$\eta(z)$ is  holomorphic and does not vanish on $\mathbb{H}$.

 A function $f(z)$ is called eta-quotient if it can
be written in the form of
\[
f(z)=\prod_{\delta |N}\eta^{r_{\delta}}(\delta z),
\]
where $N\geq 1$ and each $r_{\delta}$ is an integer. The following
two facts  is useful  to verify  whether an eta-quotient is a
modular form, see Ono \cite[p.18]{Ono04}.

\begin{prop}\label{prop2.1} If $f(z)=\prod_{\delta|N}\eta^{r_{\delta}}(\delta
z)$ is an eta-quotient with
\[ k=\frac{1}{2}\sum_{\delta|N}r_{\delta}\in \mathbb{Z},\] satisfies
the following conditions:
\begin{equation}\label{con1}
\sum_{\delta|N}\delta r_{\delta}\equiv 0 \ ({\rm mod}\ 24)
\end{equation}
and
\begin{equation}\label{con2}
\sum_{\delta|N}\frac{N}{\delta} r_{\delta}\equiv 0 \ ({\rm mod}\
24),
\end{equation}
then $f(z)$ satisfies
\begin{equation}\label{relation1}
f\left(\frac{az+b}{cz+d}\right)=\chi(d)(cz+d)^kf(z)
\end{equation}
for each $\left(
           \begin{array}{cc}
             a & b \\
             c & d \\
           \end{array}
         \right)\in \Gamma_0(N)$.
          Here the character $\chi$ is
         defined by $\chi(d):=\left(\frac{(-1)^ks}{d}\right)$, where
         \[ s:=\prod_{\delta|N}\delta^{r_{\delta}}\] and
         $\left(\frac{m}{n}\right)$ is Kronecker symbol.
\end{prop}

Based on this proposition, for a given eta-quotient $f(z)$, by
checking the conditions \eqref{con1} and \eqref{con2}, one can show
that $f(z)$ satisfies \eqref{relation1}. Moreover, if $k$ is a
positive integer and $f(z)$ is holomorphic at the cusps of
$\Gamma_0(N)$, then $f(z)\in M_k(\Gamma_0(N),\chi)$ because
$\eta(z)$ is holomorphic and does not vanish on $\mathbb{H}$.
Combined with the following proposition which gives the analytic
orders of an eta-quotient at the cusps of $\Gamma_0(N)$, we can
deduce that  $f(z)$ is a modular form.

\begin{prop}\label{prop2.2}
Let $c,d$ and $N$ be positive integers with $d|N$ and $(c,d)=1$. If
$f(z)$ is an eta-quotient satisfying the conditions in Proposition
\ref{prop2.1} for $N$, then the order of vanishing of $f(z)$ at the
cusp $\frac{c}{d}$ is
\[
\frac{N}{24}\sum_{\delta
|N}\frac{(d,\delta)^2r_{\delta}}{(d,\frac{N}{d})d\delta}.
\]
\end{prop}

In the other words, to prove that the above function $f(z)$ is holomorphic at the
cusp $\frac{c}{d}$, it suffices to check that
\[ \sum_{\delta |N}\frac{(d,\delta)^2r_{\delta}}{\delta}\geq 0.\]

Let $M$ be a positive integer and \[f(z)=\sum\limits_{n=0}^\infty
a(n)q^n\] be a function with rational integer coefficients. Define
$ord_M(f(z))$  to be the smallest $n$ such that $a(n) \not \equiv 0 \
({\rm mod\ }M)$. Sturm \cite{Sturm} provided the following powerful
criterion to determine whether two modular forms are congruent
modulo a prime by the verification of a finite number of cases.

\begin{prop}\label{prop2.3}
Let p be a prime and  $f(z), g(z)\in M_k(\Gamma_0(N))$ with rational
integer coefficients. If
\[
ord_p(f(z)-g(z))> \frac{kN}{12}\prod_{d}(1+\frac{1}{d}),
\]
 where the product is over
the prime divisors $d$ of $N$. Then $f(z)\equiv g(z)\ ({\rm mod\
}p)$, i.e., $ord_p(f(z)-g(z))=\infty$.
\end{prop}

We also need the following  result due to Lovejoy \cite{lovejoy}.

\begin{prop}\label{lem2.1}
Let \[ f=\sum\limits_{n=0}^\infty u(n)q^n\] and
\[ g=1+\sum\limits_{n=1}^\infty v(mn)q^{mn}.\] Define $w(n)$ by
\[ fg=\sum\limits_{n=0}^\infty w(n)q^n.\]  Let $d$ be a residue class
modulo $m$. Then,
\begin{description}
  \item(1) If $u(mn+d)\equiv 0\ ({\rm mod\ M})$ for $0\leq n \leq N$,
  then $w(mn+d)\equiv 0\ ({\rm mod\ M})$ for $0\leq n \leq N$.
  \item(2) If $w(mn+d)\equiv 0\ ({\rm mod\ M})$ for all $n$, then $u(mn+d)\equiv 0\ ({\rm mod\ M})$
  for all $n$.
\end{description}
\end{prop}

The following two propositions will also be used  to construct modular
forms, see Koblitz \cite{Koblitz}.

\begin{prop}\label{prop2.4}
Suppose $f(z)\in M_k(\Gamma_0(N))$ with Fourier expansion \[
f(z)=\sum\limits_{n=0}^\infty u(n)q^n.\] Then
  for any positive integer $m|N$,\[f(z)|U(m):=\sum_{n=0}^\infty u(mn)q^n\]
  is the Fourier expansion of a modular form in $M_k(\Gamma_0(N))$.
\end{prop}

\begin{prop}\label{prop2.5}
Let $\chi_1$ be a Dirichlet character modulo $M$, and let $\chi_2$
be a primitive Dirichlet character modulo $N$. Let
\[ f(z)=\sum\limits_{n=0}^\infty a(n)q^n\in M_k(M,\chi_1)\]  and
\[ g(z)=\sum\limits_{n=0}^\infty a(n)\chi_2(n)q^n. \] Then $g(z)\in
M_k(MN^2,\chi_1\chi_2^2)$. In particular, if $f(z)\in
M_k(\Gamma_0(M))$ and $\chi_2$ is quadratic, then $g(z)\in
M_k(\Gamma_0(MN^2))$.
\end{prop}

\section{Congruences for the Number of Cubic Partitions}

In this section, we give the proofs  of  Theorem \ref{thm1.3} and
Theorem \ref{thm1.4} using the technique of modular forms due to Ono
\cite{Ono96}. The following congruence relation is well-known, see,
for example, Ono \cite{Ono96}. We include a proof for the sake of
completeness.

\begin{lem}\label{lem3.1}
If $p\geq 3$ is a prime, then
\begin{equation}\label{congforp}
\frac{(q;q)_\infty ^p}{(q^p;q^p)_\infty}\equiv 1\ ({\rm mod\ }p).
\end{equation}
\end{lem}
\pf Using the well-known binomial theorem
\[
(1+x)^n=\sum_{k=0}^n {n\choose k} x^k,
\]
 it is easily seen  that
\[
\frac{(1-X)^p}{1-X^p}=\frac{\sum_{k=0}^p {p\choose
k}(-1)^kX^k}{1-X^p} \equiv \frac{1-X^p}{1-X^p} \equiv 1\ ({\rm mod\
}p),
\]
since $p\,| {p\choose k}$ for $0<k<p$. It follows that
\[
\frac{(q;q)_\infty ^p}{(q^p;q^p)_\infty}=\prod_{k=1}^\infty\frac{
(1-q^k)^p}{(1-q^{kp})} \equiv 1 \ ({\rm mod\ }p),
\]
as desired. \qed

We first consider Theorem \ref{thm1.3}, that is, for $n\geq 0$,
\[
a(25n+22)\equiv 0\ ({\rm mod}\ 5).
\]

\noindent {\it Proof.} To establish the claimed
 congruence relation, we shall construct an eta-quotient with the following expansion in $q=e^{2\pi i z}$,
\[ g(z)=\sum_{n\geq 0}  b(n)q^n.\]  We assume that $g(z)$
satisfies the following conditions
\begin{enumerate}
\item[\rm{(1)}] $g(z)$ is a modular form;
  \item[\rm{(2)}] If for all $n\geq 0$,
$ b(25n+25)\equiv 0 \ ({\rm mod}\ 5)$ then $ a(25n+22)\equiv 0 \ ({\rm mod}\ 5)$;
  \item[\rm{(3)}] The function \[ g(z)|U(25)=\sum_{n\geq 0} b(25n)q^n \equiv 0\ ({\rm mod}\ 5)\]  is also a modular
  form.
\end{enumerate}

In order to satisfy the second condition, we consider the function $h(q)$ of the following form
\[
h(q)=\prod_{i}(q^{25r_i};q^{25r_i})^{s_i}_\infty \prod_j
\Big(\frac{(q;q)_\infty^5} {(q^5;q^5)_\infty}\Big)^{t_j},
\]
where $r_i, s_i,t_j$  are integers. By the above Lemma \ref{lem3.1},
it is easily seen that for any integers $r_i$, $s_i$ and $t_j$ the
expansion of $h(q)$  has the following form modulo $5$,
\[ h(q) \equiv
1+\sum_{m\geq 1}c(m) q^{25m}\ ({\rm mod}\ 5),\] where $c(m)$ are
integers.  Now, we set
\begin{equation} \label{gzh}
g(z)=h(q) \sum_{n=0}^\infty a(n)q^{n+3}=h(q) \sum_{n\geq 3}^\infty
a(n-3)q^{n}.
\end{equation}
Since $h(q)$ is a series
in $q^{25}$ modulo $5$ with constant term $1$, we can make use of Proposition \ref{lem2.1} (2) to
 deduce that  if for any $n$, $b(25n+25)\equiv 0 \ ({\rm mod}\ 5)$, then we have $a(25n+22)\equiv 0 \ ({\rm
mod}\ 5)$.

We now proceed to determine the parameters $r_i, s_i$ and $t_j$ in $g(z)$ to make it a modular form.
 Consider the case $r_1=1, s_1=1, r_2=2, s_2=1$ and $t_1=2$, namely,
 \begin{equation}\label{gz}
 g(z)=(q^{25};q^{25})_\infty(q^{50};q^{50})_\infty \Big(\frac{(q;q)_\infty^5}
 {(q^5;q^5)_\infty}\Big)^2 \sum_{n=0}^\infty a(n)q^{n+3}.
 \end{equation}
 We are going to show that $g(z)$ satisfies the
conditions \eqref{con1} and \eqref{con2} in Proposition
\ref{prop2.1}.
Recalling
the definition of  $\eta(z)$, we can rewrite $g(z)$ as an
eta-quotient
\begin{eqnarray*}
g(z) &=&\frac{\eta(25z)\eta(50z)}{\eta(z)\eta(2z)}
\left(\frac{\eta^5(z)}{\eta(5z)}\right)^2 \\[5pt]
&=& \frac{\eta^9(z)\eta(25z)\eta(50z)}{\eta(2z)\eta^2(5z)}.
\end{eqnarray*}
The two conditions \eqref{con1} and \eqref{con2} can be expressed as follows,
\begin{eqnarray*}
\sum_{\delta | 50} \delta r_{\delta}&=&9-2-5\times 2+25+50 \, \equiv  0\ (\rm{mod}\ 24),\\[5pt]
\sum_{\delta | 50} \frac{50}{\delta} r_{\delta}&=&50\times
9-25-10\times 2+2+1\equiv 0\ (\rm{mod}\ 24).
\end{eqnarray*}

To make $g(z)$ a modular form, it remains to
compute the order of $g(z)$ at cusps. By
Proposition \ref{prop2.2}, it is easily verified  that the order
of $g(z)$ at the cusps of
 $\Gamma_0(50)$ are nonnegative, that is, for any $d|50$,
 \[
 \sum_{\delta| 50} \frac{(d,\delta)^2 r_\delta}{\delta} \geq 0.\]
So we have
$g(z)\in M_4(\Gamma_0(50),\chi_1)$, where
\[ \chi_1(d)=\left(\frac{25}{d}\right)=\left(\frac{5}{d}\right)^2=1\]
for $(d,50)=1$. This implies that $g(z)\in M_4(\Gamma_0(50))$.

Since we have proved that $g(z)$ satisfies the second condition,
 the following congruence is valid
\begin{equation}\label{conga5}
 a(25n+22)\equiv 0\ ({\rm mod\ 5})
\end{equation}
 provided that for all $n \geq 0$,
\begin{equation}\label{congb5}
b(25n+25)\equiv 0\ ({\rm mod\ 5}).
\end{equation}
Let us rewrite (\ref{congb5}) as
\begin{equation}\label{congbq}
\sum_{n\geq 1}b(25n)q^n \equiv 0\ ({\rm mod\ 5}).
\end{equation}
By Proposition \ref{prop2.4}, we see that the summation on
the left hand side of (\ref{congbq}) is a modular form, that is,
\[
g(z)|U(25)=\sum_{n\geq 1}b(25n)q^n \in M_4(\Gamma_0(50)).
\]
Hence, by Proposition \ref{prop2.3}, we find that \eqref{congbq} is
valid if \eqref{congb5} holds for \[0\leq n\leq \frac{4\times
50}{12}(1+\frac{1}{2})(1+\frac{1}{5})+1=31.\]

Applying Lemma \ref{lem3.1} with $p=5$ to \eqref{gz}, we have
\begin{equation}\label{bconga}
\sum_{n\geq 3}b(n)q^n \equiv
(q^{25};q^{25})_\infty(q^{50};q^{50})_\infty \sum_{n=0}^\infty
a(n)q^{n+3}\ ({\rm mod}\ 5).
\end{equation}
Using the above relation and Proposition \ref{lem2.1} (1),
 we see that the verification of \eqref{congb5} on $b(n)$ for $0\leq n \leq 31$
 can be reduced to the verification of \eqref{conga5} on $a(n)$ for
a the same range $0 \leq n \leq 31$.
 It is
readily checked that \eqref{conga5} holds for $0\leq n \leq 31$.
This completes the proof. \qed

Now, we turn to the proof of  Theorem \ref{thm1.4}, namely,
\begin{equation}\label{4c}
a(49n+15)\equiv a(49n+29) \equiv a(49n+36) \equiv a(49n+43) \equiv 0
\ ({\rm mod}\ 7).
\end{equation}

\noindent {\it Proof of Theorem \ref{thm1.4}.} Following the above
procedure in the proof of Theorem \ref{thm1.3}, by the generating
function of $a(n)$, we  construct an eta-quotient
\begin{eqnarray}
h(z)&=&\frac{\eta(49z)\eta(98z)}{\eta(z)\eta(2z)}
\left(\frac{\eta^7(z)}{\eta(7z)}\right)^2\nonumber\\[5pt]
&=&\frac{\eta^{13}(z)\eta(49z)\eta(98z)}{\eta(2z)\eta^2(7z)}\nonumber\\[5pt]
&=& \left(\frac{(q;q)_\infty
^7}{(q^7;q^7)_\infty}\right)^2(q^{49};q^{49})_\infty(q^{98};q^{98})_\infty
\sum_{n=0}^\infty a(n)q^{n+6}.\label{hz}
\end{eqnarray}

Setting $N=98$,  we see that $h(z)$ satisfies the conditions
\eqref{con1} and \eqref{con2} in Proposition \ref{prop2.1}, namely,
\begin{eqnarray*}
\sum_{\delta | 98} \delta r_{\delta}&=&13-2-2\times 7+49+98\, \equiv  0\ (\rm{mod}\ 24),\\[5pt]
\sum_{\delta | 98} \frac{98}{\delta} r_{\delta}&=&98\times
13-49-14\times 2+2+1\equiv 0\ (\rm{mod}\ 24).
\end{eqnarray*}  Moreover, it is not difficult to verify that the order
of $h(z)$ at the cusps of
 $\Gamma_0(98)$ are nonnegative by Proposition \ref{prop2.2}, that is, for any $d|98$,
 \[
 \sum_{\delta| 98} \frac{(d,\delta)^2 r_\delta}{\delta} \geq 0.
 \]
 Hence we deduce that $h(z)$ is a modular form, i.e.,
$h(z)\in M_6(\Gamma_0(98),\chi_2)$. Moreover,
\[ \chi_2(d)=\left(\frac{49}{d}\right)=\left(\frac{7}{d}\right)^2=1\]
for $(d,98)=1$. This implies that  $h(z)\in M_6(\Gamma_0(98))$.

Write
\[
h(z)=\sum_{n\geq 6}c(n)q^n.
\]
Applying Lemma \ref{lem3.1} with $p=7$, \eqref{hz} becomes
\begin{equation}\label{cconga}
\sum_{n\geq 6}c(n)q^n \equiv
(q^{49};q^{49})_\infty(q^{98};q^{98})_\infty \sum_{n\geq 6}
a(n-6)q^{n}\ ({\rm mod}\ 7).
\end{equation}
Since $(q^{49};q^{49})_\infty(q^{98};q^{98})_\infty$ can be expanded as series
in $q^{49}$ with constant term $1$, we can make use of   Proposition \ref{lem2.1} (2) to
 deduce that the four congruences in (\ref{4c}) can be derived from
 the corresponding congruences for $c(n)$, i.e., for $n \geq 0$,
\begin{equation}\label{congc}
c(49n+21)\equiv c(49n+35) \equiv c(49n+42) \equiv c(49n+49) \equiv 0
\ ({\rm mod}\ 7).
\end{equation}

Observing that the arithmetic progressions $49n+21, 49n+35, 49n+42, 49n+49$
in the above congruences  are divided by $7$, we construct another function based on
 $c(n)$ as follows
\[ u(z)=\sum_{n\geq 1}d(n)q^n=\sum_{n\geq 1}c(7n)q^n.
\]
By Proposition \ref{prop2.4}, we have $ u(z) \in M_6(\Gamma_0(98))$.
Obviously,  \eqref{congc} can be restated as
\begin{equation}\label{congd}
d(7n+3)\equiv d(7n+5) \equiv d(7n+6) \equiv d(7n+7) \equiv 0 \ ({\rm
mod}\ 7).
\end{equation}

As will be seen, one can combine the above four congruences into a single
congruence relation.
Define
\begin{equation}
 v(z)=\sum_{n\geq 1} e(n)q^n=\sum_{n\geq 1}d(n)q^n-
\sum_{n\geq 1} \left(\frac{n}{7}\right)d(n)q^n \label{k1z1}.
\end{equation}
Since $\left(\frac{n}{7}\right)=1$ for  $n\equiv 1,2,4\ ({\rm
mod}\ 7)$, $\left(\frac{n}{7}\right)=-1$ for $n\equiv 3,5,6\ ({\rm
mod}\ 7)$ and  $\left(\frac{n}{7}\right)=0$ for $n\equiv 0\ ({\rm
mod}\ 7)$,  $v(z)$ can be expressed as
\begin{equation}
v(z)=\sum_{\left(\frac{n}{7}\right)=-1}2d(n)q^n+\sum_{n\equiv 0\ ({\rm
mod}\ 7)}d(n)q^n.\label{k1z}
\end{equation}
To prove
\eqref{congd}, it suffices to show that $v(z) \equiv 0 \ ({\rm mod}\
7)$.

Denote the second summation in \eqref{k1z1} by \[
w(z)=\sum_{n\geq 1} \left(\frac{n}{7}\right)d(n)q^n. \]  By Proposition \ref{prop2.5},
we see that $w(z)$ is a modular form. In other words, $
w(z)\in M_6(\Gamma_0(4802))$. Since $u(z)$ is also a modular form, we obtain
\[ v(z)=u(z)-w(z)\in M_6(\Gamma_0(4802).\]
By Proposition
\ref{prop2.3}, we find that $v(z)\equiv 0\ ({\rm mod\ }7)$ can be verified by a finite number of cases.
To be precise, we need to check that
\[
e(n)\equiv 0\ ({\rm mod\ 7})
\]
holds for $0\leq n\leq 4117$.  In view of \eqref{k1z}, we only need to verify that
 \eqref{congd} holds for $0 \leq n \leq
\lceil\frac{4117-1}{7}\rceil-1=587$. Since $d(n)=c(7n)$,
it suffices to check \eqref{congc}  for $0 \leq n \leq 587$.
Finally, using \eqref{cconga} and Proposition \ref{lem2.1} (1),
it is necessary to verify  Theorem \ref{thm1.4} holds only for $0\leq
n\leq 587$, which is an easy task. This
completes the proof. \qed

\section{The Parity of $a(n)$}

In this section, we show that the function $a(n)$ takes infinitely
many even values and infinitely many odd values.  For the partition
function $p(n)$, it has been conjectured by Parkin and Shanks
\cite{Parkin67} that the parities of  $p(1), p(2), \ldots, p(N)$ are
equidistributed when $N$ tends to infinity. Using the Euler's
recurrence formula for $p(n)$,
\begin{equation}\label{euler}
p(n)+\sum_{0<\omega_j\leq n}(-1)^{j}p(n-\omega(j))=0,
\end{equation}
where $\omega(j)=j(3j-1)/2,-\infty<j<\infty$, Kolberg \cite{Kolberg}
proved that $p(n)$ takes both even and odd values infinitely often.

To prove the analogous property for $a(n)$, we need Jacobi's identity
\begin{equation}\label{jacobi1}
\sum_{n=0}^\infty (-1)^n(2n+1)q^{n(n+1)}=(q^2;q^2)_\infty ^3
\end{equation}
and Gauss's identity
\begin{equation}\label{gauss1}
\sum_{n=0}^\infty q^{n(n+1)/2}=\frac{(q^2;q^2)_\infty
^2}{(q;q)_\infty}.
\end{equation}
We now have the following recurrence relation modulo 2.

\begin{thm}\label{thm4.1}
\begin{equation}\label{recurrence1}
a(n)+\sum_{0<k+k^2\leq n}a(n-k-k^2)  \equiv \Delta(n) \ ({\rm mod}\
2),
\end{equation}
where $\Delta(n)=1$, if $n=s(s+1)/2$ for some integer $s$ and
$\Delta(n)=0$, otherwise.
\end{thm}

\pf Multiplying both sides of
\eqref{Defa} by $(q^2;q^2)_\infty^3$, we get
\begin{equation}\label{a(n)}
(q^2;q^2)_\infty ^3\sum_{n=0}^\infty
a(n)q^n=\frac{(q^2;q^2)_\infty^2}{(q;q)_\infty}.
\end{equation}
Substituting \eqref{jacobi1} and \eqref{gauss1} into
both sides of \eqref{a(n)}, we find that
\begin{eqnarray}
\sum_{n=0}^\infty
q^{n(n+1)/2}&=&\sum_{n=0}^\infty(-1)^n(2n+1)q^{n(n+1)}\sum_{n=0}^\infty
a(n)q^n\nonumber \\[5pt]
&\equiv& \sum_{n=0}^\infty q^{n(n+1)}\sum_{n=0}^\infty a(n)q^n \quad
({\rm mod}\ 2)\label{conga(n)}.
\end{eqnarray}
Equating coefficients of $q^n$ on both sides of
 \eqref{conga(n)} gives   \eqref{recurrence1}. This
completes the proof. \qed

With the aid of the formula \eqref{recurrence1} and following
the idea of Kolberg \cite{Kolberg} for $p(n)$, we obtain the following theorem for
$a(n)$.

\begin{thm}\label{thm1.6}
There are infinitely many integers $n$ such that  $a(n)$ is even
and there are infinitely many integers $n$ such that  $a(n)$ is odd.
\end{thm}

\pf We prove by contradiction. Assume that there exists $m$ such
that $a(n)$ is odd for any $n\geq m$.   Without loss of generality,
we may assume that  $m$ is an even integer greater than $2$. It is
easy to show that there exists an integer $m^2+2m\leq t \leq
m^2+3m+1$
 such that $\Delta(t)=0$. Setting
$t=m^2+2m+\delta$, where $0\leq \delta \leq m+1$. Substituting
$n=t$ into \eqref{recurrence1} yields
 \begin{equation}\label{parity}
 a(m^2+2m+\delta)+
a(m^2+2m+\delta-2)+\cdots+a(m+\delta)\equiv 0 \ ({\rm mod}\ 2).
\end{equation}
But the left hand side of \eqref{parity} is the sum of  $m+1$ odd
numbers,  so it is also odd since $m$ is even. This leads to a
contradiction with parity of the right hand side of \eqref{parity}.

On the other hand, assume that $a(n)$ is even for any $n\geq m$,
where $m\geq 5$. It is easy to verify the following inequalities for
$k\geq 10$,
\[
k(k-1)<\frac{\sqrt{2}k\big(\sqrt{2}k-1\big)}{2}<\frac{\sqrt{2}k\big(\sqrt{2}k+1\big)}{2}<k(k+1).
\]
Therefore, for $k\geq 10$, there are no integers $e(e+1)$ in the following interval
with $\lceil\sqrt{2}k\rceil$ elements,
\[
\left[\frac{\lceil\sqrt{2}k\rceil\big(\lceil\sqrt{2}k\rceil-1\big)}{2},\frac{\lceil\sqrt{2}k\rceil\big(\lceil\sqrt{2}k\rceil+1\big)}{2}\right].
\]
Choose $k$ such that $d=\lceil\sqrt{2}k\rceil>m$ and set \[
t=\frac{d(d+1)}{2}.\]  Let $r$ be the largest integer $k$
such that $k^2+k\leq t$. It follows  that
\[t-r-r^2 > d >m.\] Substituting $n=t$ into
\eqref{recurrence1}, we obtain that
\[
a(t)+a(t-2)+\cdots+a(t-r-r^2)\equiv
\Delta(t)\equiv 1\ ({\rm mod}\ 2).
\]
This is impossible since the left hand side of above
congruence are a sum of even numbers. This completes the proof. \qed

\vspace{.2cm} \noindent{\bf Acknowledgments.} We wish to thank Lisa H. Sun for helpful comments.
This work was
supported by the 973 Project, the PCSIRT Project of the Ministry of
Education, and the National
Science Foundation of China.

\end{document}